\documentclass[journal]{IEEEtran}
\usepackage{amssymb}
\usepackage{epsfig}
\usepackage{graphicx}
\usepackage{colortbl}
\usepackage{float}
\usepackage{multirow}

\newtheorem{theorem}{Theorem}

\newtheorem{remark}{Remark}
\newtheorem{lemma}{Lemma}
\newtheorem{assumption}{Assumption}

\usepackage{algorithm}
\usepackage{algorithmic}
\ifCLASSINFOpdf
\else
\fi
\begin{document}
    %
    \title{Decentralized Control of Linear Systems with Private Input and Measurement Information}
    %
    %
    %

    \author{Juanjuan~Xu and~Huanshui~Zhang 
        \thanks{*This work was supported by the National Natural Science Foundation of China under Grants 61821004, 62250056
            and the Natural Science Foundation of Shandong Province
            (ZR2021ZD14, ZR2021JQ24, ZR2020ZD24), and High-level
            Talent Team Project of Qingdao West Coast New Area
            (RCTD-JC-2019-05), Key Research and Development Program of Shandong Province (2020CXGC01208), and Science
            and Technology Project of Qingdao West Coast New Area
            (2019-32, 2020-20, 2020-1-4).}
        \thanks{J. Xu is with School of Control Science and Engineering, Shandong University, Jinan, Shandong, P.R. China 250061.
                {\tt\small juanjuanxu@sdu.edu.cn}}
       \thanks{H. Zhang is with the College of Electrical Engineering and Automation, Shandong University of Science and Technology, Qingdao, Shandong, P.R. China 266000.
            {\tt\small hszhang@sdu.edu.cn}}
    }

    \maketitle

    \begin{abstract}

        In this paper, we study the linear quadratic (LQ) optimal control problem of linear systems with private input and measurement information.
        The main challenging lies in the unavailability of other regulators' historical input information. To overcome this difficulty, we introduce
        a kind of novel observers by using the private input and measurement information and accordingly design a kind of new decentralized controllers.
        In particular, it is verified that the corresponding cost function under the proposed decentralized controllers are asymptotically optimal as comparison with the optimal cost under optimal state-feedback controller.
 The presented results in this paper are new to the best of our knowledge, which represent the fundamental contribution to classical decentralized control.

    \end{abstract}


    %
    \IEEEpeerreviewmaketitle

    \section{Introduction}

    In the past decades, control problems of large-scale systems have
    received increasing attention due to the wide applications in practical fields including
    chemical engineering system, aerospace system, power industry and
    telecommunication network \cite{Largescale1}-\cite{Largescale2}. In such systems,
    there exist a number of independent subsystems which serve particular functions, share resources,
    and are governed by a set of interrelated goals and constraints \cite{Largescale}. Generally,
    there are two kinds of control methods. One is the centralized which requires a central station and a communication network
    to transfer information between subsystems and the central station. Considering the physical
    constraints like high cost of cabling and limited communication bandwidth, centralized controllers are difficult to implement
    and therefore not desirable. The other is decentralized which only applies the available local input and output information
    and is thus popular.

    Plenty of progress has been made in decentralized control on the stabilization and linear quadratic Gaussian (LQG) problems \cite{Literature1}-\cite{Literature5}. Firstly, for the decentralized stabilization of the large-scale systems, algebraic necessary and sufficient
    conditions have been obtained for making a singular system both regular and impulse free by decentralized output feedback control laws and decentralized proportional-plus-derivative output feedback control in \cite{Output}. \cite{Zak} studied the problem of robust stabilization for a class of nonlinear interconnected systems consisting of linear subsystems coupled by nonlinear interconnections that are unknown and quadratically
    bounded by using a decentralized dynamic output feedback based linear controller.
    In \cite{TongShaocheng}, the problem of adaptive fuzzy decentralized output-feedback control design was investigated for a class of switched nonlinear large-scale systems in strict-feedback form.
    \cite{Miller1} obtained that if the system is centrally
    controllable and observable and the graph associated with the system
    is strongly connected, then there exists a linear periodic controller
    which provide LQR performance as close as desired to the optimal centralized performance.
    Furthermore, \cite{Miller} considered the problem of LQR optimal control in the context of decentralized systems where the system matrices are parameterized by a constant larger than one, and demonstrated that the optimal decentralized performance can be arbitrarily bad in comparison to the optimal centralized performance even when the plant has no unstable
    decentralized fixed modes.
    Secondly, for the decentralized LQG problem, \cite{Onestep} obtained the optimal solutions to the case with
    one-step delayed information sharing pattern by using dynamic programming. In \cite{Onestep1}, the Nash equilibrium was given for the multistage LQG decision problems with two decision makers and one-step delay observation sharing
    pattern. A partial historical sharing was studied in \cite{Nayyar1}, and \cite{Nayyar} gave the optimal linear
    control strategies by applying the common information approach.
    The optimal linear decentralized control was given for the $d$-step delayed information sharing pattern in \cite{dStep}.
    It is noted that the controllers in the above mentioned work share their all or partial historical control information where the estimation/predictor is least square one based on the local information.
    However, when the control information is unavailable to the other regulators, how to design the observer or estimation remains to be solved.

    In this paper, we study the linear quadratic (LQ) optimal control problem of linear systems with private input and measurement information, that is,
    both the measurement and the control information are unavailable to the other regulators which makes the solvability much more difficult.
    To this end, we define a kind of novel observers by using the private input and measurement information
    and then design a kind of new decentralized controllers. It is further verified that the corresponding cost function under the proposed decentralized controllers are asymptotically optimal as comparison with the optimal cost under optimal state-feedback controller.

    The remainder of the paper is organized as follows. Section II presents the studied problem.
    Some preliminaries on state-feedback control is given in Section III. The decentralized control is shown in Section IV for the problem with input sharing pattern. Section V presents the decentralized control with private input information.
 The results for the system with multiple inputs are shown in Section VI.
    Numerical examples are given in Section VII. Some concluding remarks are given in the last section.

    The following notations will be used throughout this paper: $R^n$
    denotes the set of $n$-dimensional vectors; $x'$ denotes the
    transpose of $x$;
    a symmetric matrix $M>0\ (\geq 0)$ means that $M$ is
    strictly positive-definite (positive semi-definite). $I$ and $0$ denotes the identity matrix and the zero vector/matrix
    with appropriate dimensions.
  $\|x\|\triangleq \sqrt{x'x}$ means the 2-norm of vector $x$, and $\|M\|\triangleq \sqrt{\lambda_{max}(M'M)}$ means the induced 2-norm of matrix $M$.


    \section{Problem Formulation}

    In this paper, we consider the following linear system:
    \begin{eqnarray}
        x(k+1)&=&Ax(k)+B_1u_1(k)+B_2u_2(k),\label{o1}\\
        y_i(k)&=&H_ix(k), i=1,2,\label{o2}
    \end{eqnarray}
    where $x\in R^n$ is the state, $y_i\in R^{s_i}$ is the measurement, $u_i\in R^{m_i}$ is the control input with $s_i, m_i$ being positive integers.
    $A, B_i$ and $H_i$ are time-invariant matrices with compatible dimensions.
    The initial value is given by $x(0)=x_0$.

    The cost function is given by
    \begin{eqnarray}
        J&=&\sum_{k=0}^\infty[x'(k)Qx(k)+u_1'(k)R_1u_1(k)\nonumber\\
        &&+u_2'(k)R_2u_2(k)],\label{o5}
    \end{eqnarray}
    where $Q, R_i$ are positive semi-definite matrices with compatible dimensions.

    Considering that two controllers $u_1, u_2$ and two measurements $y_1, y_2$ involved in the system (\ref{o1}),
    there exist several kinds of information patterns to design the controllers such as state-feedback, decentralized with control information sharing pattern, and decentralized with private input information.
    Recalling that the main difficulty lies in the unavailability of the other regulators' control information, we aim to design a kind of observers and decentralized controllers by using the private input and measurement information in the sequel. Mathematically, by denoting
    \begin{eqnarray}
        Y_i(k)=\{y_i(0), \ldots, y_i(k), u_i(0), \ldots, u_i(k-1)\},\label{o3}
    \end{eqnarray}
    we will design a kind of novel observer-feedback decentralized controllers $u_i(k)$ based on the information $Y_i(k)$ for $i=1,2$.

    To this end, we make the following assumption:
    \begin{assumption}\label{a1}
        System $(A,B)$ is stabilizable and system $(A,Q)$ is observable where $B=\left[
        \begin{array}{cc}
            B_1 & B_2 \\
        \end{array}
        \right].
        $
    \end{assumption}

    \section{Preliminaries on State-Feedback Control}\label{sec1}

    The most ideal scenario is that both the state information $\{x(0), \ldots, x(k)\}$ and the input information $\{u_i(0), \ldots, u_i(k-1), i=1,2\}$ are available to design the controllers $u_1(k)$ and $u_2(k)$.
    In this case, the solvability of the optimal control problem of minimizing the cost function (\ref{o5}) subject to (\ref{o1}) is reduced to a standard LQ optimal control problem \cite{Anderson}.
    \begin{lemma}
        Under Assumption \ref{a1}, the optimal solution of minimizing (\ref{o5}) subject to (\ref{o1}) is given by
        \begin{eqnarray}
            u^*(k)&=&Kx(k),\label{o6}
        \end{eqnarray}
        where the feedback gain $K$ is given by
        \begin{eqnarray}
            K&=&-(R+B'PB)^{-1}B'PA,\label{o7}
        \end{eqnarray}
        and $P$ satisfies the algebraic Riccati equation (ARE):
        \begin{eqnarray}
            P&=&A'PA+Q-A'PB(R+B'PB)^{-1}B'PA,\label{o8}
        \end{eqnarray}
        with $R=\left[
        \begin{array}{cc}
            R_1 & 0 \\
            0 & R_2 \\
        \end{array}
        \right].$

        The optimal cost function is
        \begin{eqnarray}
            J^*=x'(0)Px(0).\label{o22}
        \end{eqnarray}
    \end{lemma}


    \section{Decentralized Control with Input Sharing Pattern}\label{sec2}

    The second scenario is that the measurement information is private and the input information is sharing, that is,
    by denoting
    \begin{eqnarray}
        \bar{Y}_i(k)&=&\{y_i(0), \ldots, y_i(k), u_1(0), \ldots, u_1(k-1), \nonumber\\
        &&~~~~~~~~~~~~~~u_2(0), \ldots, u_2(k-1)\}, \label{o25}
    \end{eqnarray}
    the controller $u_i(k)$ is $\bar{Y}_i(k)$-causal. Considering the sharing of the input information,
    we design the observers $\hat{x}_i^\diamond, i=1,2$ as below:
    \begin{eqnarray}
        \hat{x}_1^\diamond(k+1)&=&A\hat{x}_1^\diamond(k)+B_1u_1^\diamond(k)+B_2u_2^\diamond(k)\nonumber\\
        &&+L_1^\diamond[y_1(k)-H_1\hat{x}_1^\diamond(k)],\label{o29}\\
        \hat{x}_2^\diamond(k+1)&=&A\hat{x}_2^\diamond(k)+B_1u_1^\diamond(k)+B_2u_2^\diamond(k)\nonumber\\
        &&+L_2^\diamond[y_2(k)-H_2\hat{x}_2^\diamond(k)],\label{o30}
    \end{eqnarray}
    where the observer gains $L_1^\diamond$ and $L_2^\diamond$ are to be chosen such that the observers are stable.

    \begin{theorem}\label{t3}
        Assume that system $(A, H_i)$ is observable for $i=1, 2$, then there exist observer gains $L_i^\diamond, i=1, 2$
        such that
        \begin{eqnarray}
            \lim_{k\rightarrow\infty}\tilde{x}^\diamond(k)=0,\label{o28}
        \end{eqnarray}
        where
        \begin{eqnarray}
            \tilde{x}^\diamond(k)=\left[
            \begin{array}{c}
                \tilde{x}_1^\diamond(k) \\
                \tilde{x}_2^\diamond(k) \\
            \end{array}
            \right]\triangleq  \left[
            \begin{array}{c}
                x(k)-\hat{x}_1^\diamond(k) \\
                x(k)-\hat{x}_2^\diamond(k) \\
            \end{array}
            \right].
        \end{eqnarray}
    \end{theorem}
    \emph{Proof.} From (\ref{o1}), (\ref{o29}) and (\ref{o30}), it follows that
    \begin{eqnarray}
        \tilde{x}^\diamond(k+1)=\check{A}\tilde{x}^\diamond(k).\label{o32}
    \end{eqnarray}
    where $\check{A}=\left[
    \begin{array}{cc}
        A-L_1H_1 & 0 \\
        0 & A-L_2H_2 \\
    \end{array}
    \right]$.
    From the observability of $(A, H_i)$ for $i=1, 2$, then there exists $L_i^\diamond, i=1, 2$ such that
    (\ref{o28}) holds. The proof is now completed. \hfill $\blacksquare$

    In this case, we design the decentralized controllers as follows:
    \begin{eqnarray}
        u_1^\diamond(k)=K_1\hat{x}_1^\diamond(k),\label{o26}\\
        u_2^\diamond(k)=K_2\hat{x}_2^\diamond(k),\label{o27}
    \end{eqnarray}
    where $K_1=\left[
    \begin{array}{cc}
        I & 0 \\
    \end{array}
    \right]K, K_2=\left[
    \begin{array}{cc}
        0 & I \\
    \end{array}
    \right]K$ with $K$ defined by (\ref{o7}).
    The corresponding closed-loop system (\ref{o1}) under (\ref{o26})-(\ref{o27}) is stable.
    \begin{theorem}\label{t4}
        Assume that Assumption \ref{a1} holds and system $(A, H_i)$ is observable for $i=1, 2$, then the closed-loop system (\ref{o1}) under (\ref{o26})-(\ref{o27}) is stable:
        \begin{eqnarray}
            x(k+1)&=&Ax(k)+B_1K_1\hat{x}_1^\diamond(k)+B_2K_2\hat{x}_2^\diamond(k).\label{o31}
        \end{eqnarray}
    \end{theorem}
    \emph{Proof.} From (\ref{o31}), the closed-loop system is reformulated as
    \begin{eqnarray}
        x(k+1)&=&(A+BK)x(k)+\mathcal{\mathcal{B}}\tilde{x}^\diamond(k),\nonumber
    \end{eqnarray}
 where $\mathcal{B}=\left[
    \begin{array}{cc}
        -B_1K_1 & -B_2K_2 \\
    \end{array}
    \right]$.
    Combining with (\ref{o32}), one has
    \begin{eqnarray}
        \left[
        \begin{array}{c}
            x(k+1) \\
            \tilde{x}^\diamond(k+1) \\
        \end{array}
        \right]=\tilde{A}\left[
        \begin{array}{c}
            x(k) \\
            \tilde{x}^\diamond(k) \\
        \end{array}
        \right],\nonumber
    \end{eqnarray}
    where $\tilde{A}=\left[
    \begin{array}{cc}
        A+BK & \mathcal{B} \\
        0 & \check{A} \\
    \end{array}
    \right]
    $. By using the stabilizability of $(A,B)$ and the observability of $(A, Q)$, we have that $A+BK$ is stable.
    Together with the stability of $\check{A}$, the stability of the closed-loop system (\ref{o31}) follows.
    The proof is now completed. \hfill $\blacksquare$

    \section{Decentralized Control with Private Input Information}

    Different from the cases discussed in Section \ref{sec1} and \ref{sec2} where the input information is sharing among regulators,
    we further study a scenario that all information acquired by the regulator is confidential which is unavailable to the other regulators.
    Mathematically, the design of the controller $u_i(k)$ depends on the information
    \begin{eqnarray}
        Y_i(k)=\{y_i(0), \ldots, y_i(k), u_i(0), \ldots, u_i(k-1)\}.\label{o3}
    \end{eqnarray}

    \begin{remark}\label{rem1}
        The situation that no control information sharing among regulators exists widely in cyber-physical system \cite{application}. However, considering that the fact that the information $\{u_2(0), \ldots, u_2(k-1)\}$ is not available to $u_1(k)$, the observer $\hat{x}_1^\diamond(k)$ defined by (\ref{o29}) loses efficacy. Observer $\hat{x}_2^\diamond(k)$ defined by (\ref{o30}) also
        faces the same problem. Accordingly, how to design the decentralized controllers remains challenging in the case with private
        input and measurement information.

    \end{remark}

    To overcome the unavailability of the other control information, we design novel decentralized controllers as
    \begin{eqnarray}
        u_1^\star(k)=K_1\hat{x}_1(k),\label{o11}\\
        u_2^\star(k)=K_2\hat{x}_2(k),\label{o12}
    \end{eqnarray}
    where the observers $\hat{x}_i, i=1,2$ satisfy
    \begin{eqnarray}
        \hat{x}_1(k+1)&=&A\hat{x}_1(k)+B_1u_1^\star(k)+B_2K_2\hat{x}_1(k)\nonumber\\
        &&+L_1[y_1(k)-H_1\hat{x}_1(k)],\label{o9}\\
        \hat{x}_2(k+1)&=&A\hat{x}_2(k)+B_1K_1\hat{x}_2(k)+B_2u_2^\star(k)\nonumber\\
        &&+L_2[y_2(k)-H_2\hat{x}_2(k)],\label{o10}
    \end{eqnarray}
    where the observer gains $L_1$ and $L_2$ are to be chosen such that the observers are stable.

    By denoting the matrix
    $$\mathcal{A}\triangleq \left[
    \begin{array}{cc}
        A+B_2K_2-L_1H_1 & -B_2K_2 \\
        -B_1K_1 & A+B_1K_1-L_2H_2 \\
    \end{array}
    \right],$$
    we have the detailed selection of $L_i, i=1, 2$ and the stability of the observers as shown below.
    \begin{theorem}\label{t1}
        If the observer gains $L_i, i=1,2$ are chosen such that the matrix $\mathcal{A}$ is stale,
        then the observers $\hat{x}_i(k), i=1, 2$ are stable under controllers (\ref{o11})-(\ref{o12}) in the sense that
        \begin{eqnarray}
            \lim_{k\rightarrow\infty} \|\hat{x}_i(k)-x(k)\|=0,~i=1, 2.\nonumber
        \end{eqnarray}
        Moreover, under Assumption \ref{a1}, the closed-loop system (\ref{o1}) under controllers (\ref{o11})-(\ref{o12}) is stable.
    \end{theorem}
    \emph{Proof.} Firstly, from (\ref{o1}) and (\ref{o11})-(\ref{o12}), the closed-loop system is reduced to
    \begin{eqnarray}
        x(k+1)&=&Ax(k)+B_1K_1\hat{x}_1(k)+B_2K_2\hat{x}_2(k)\nonumber\\
        &=&(A+BK)x(k)-B_1K_1\tilde{x}_1(k)-B_2K_2\tilde{x}_2(k),\nonumber\\\label{o13}
    \end{eqnarray}
    where $\tilde{x}_i(k)\triangleq x(k)-\hat{x}_i(k)$.

    Together with (\ref{o9}), it yields that
    \begin{eqnarray}
        \tilde{x}_1(k+1)
        &=&(A+BK)x(k)-B_1K_1\tilde{x}_1(k)-B_2K_2\tilde{x}_2(k)\nonumber\\
        &&-(A+BK)\hat{x}_1(k)-L_1H_1\tilde{x}_1(k)\nonumber\\
        &=&(A+B_2K_2\!-\!L_1H_1)\tilde{x}_1(k)\!-\!B_2K_2\tilde{x}_2(k),\label{o14}
    \end{eqnarray}
    By using (\ref{o10}) and (\ref{o13}), we have that
    \begin{eqnarray}
        \tilde{x}_2(k+1)
        &=&(A+BK)x(k)-B_2K_2\tilde{x}_2(k)-B_1K_1\tilde{x}_1(k)\nonumber\\
        &&-(A+BK)\hat{x}_2(k)-L_2H_2\tilde{x}_2(k)\nonumber\\
        &=&(A+B_1K_1\!-\!L_2H_2)\tilde{x}_2(k)\!-\!B_1K_1\tilde{x}_1(k).\label{o15}
    \end{eqnarray}
    Combining with (\ref{o14}) and (\ref{o15}), it follows that
    \begin{eqnarray}
        \tilde{x}(k+1)=\mathcal{A}\tilde{x}(k),\label{o16}
    \end{eqnarray}
    where $\tilde{x}(k)\triangleq\left[
    \begin{array}{c}
        \tilde{x}_1(k) \\
        \tilde{x}_2(k) \\
    \end{array}
    \right]$.

    Accordingly, if the matrix $\mathcal{A}$ is stable, then $\lim_{k\rightarrow\infty}\tilde{x}(k)=0$, that is, the observer $\hat{x}_i(k)$ is stable under (\ref{o11})-(\ref{o12}).

    Furthermore, from (\ref{o13}) and (\ref{o16}), it follows that
    \begin{eqnarray}
        \left[
        \begin{array}{c}
            x(k+1) \\
            \tilde{x}(k+1) \\
        \end{array}
        \right]=\bar{A}\left[
        \begin{array}{c}
            x(k) \\
            \tilde{x}(k) \\
        \end{array}
        \right],\label{o19}
    \end{eqnarray}
    where $\bar{A}=\left[
    \begin{array}{cc}
        A+BK & \mathcal{B}
        \\
        0 & \mathcal{A} \\
    \end{array}
    \right].$
    By using the stability of
    the matrices $A+BK$ and $\mathcal{A}$, we have that the closed-loop system (\ref{o13}) is stable.
    The proof is now completed. \hfill $\blacksquare$

    Noting that the stability of the state $x(k)$ and the observers $\hat{x}_i(k), i=1, 2$ has been guaranteed by choosing appropriate matrices $L_1$ and $L_2$ as shown in Theorem \ref{t1}. We next study the asymptotical optimal property of the decentralized controllers (\ref{o11}) and (\ref{o12}). To this end, we define the following cost function:
    \begin{eqnarray}
        J(s, M)&=&\sum_{k=s}^M[x'(k)Qx(k)+u_1'(k)R_1u_1(k)\nonumber\\
        &&+u_2'(k)R_2u_2(k)].\label{o20}
    \end{eqnarray}
    and denote \begin{eqnarray}
        S_1&=&(A+BK)'P\mathcal{B}-\left[
        \begin{array}{cc}
            K_1'R_1K_1 & K_2'R_2K_2
        \end{array}
        \right],\nonumber\\
        S_2&=&\left[
        \begin{array}{cc}
            K_1'R_1K_1 & 0 \\
            0 & K_2'R_2K_2 \\
        \end{array}
        \right]+\mathcal{B}'P\mathcal{B}.\nonumber
    \end{eqnarray}
    \begin{theorem}\label{t2}
        Under Assumption \ref{a1}, the corresponding cost function (\ref{o20}) under decentralized controllers (\ref{o11}) and (\ref{o12}) with $L_i, i=1, 2$ chosen according to Theorem \ref{t1} is given by
        \begin{eqnarray}
            &&J^\star(s,\infty)\nonumber\\
            &=&x'(s)Px(s)+\sum_{k=s}^\infty\left[
            \begin{array}{c}
                x(k) \\
                \tilde{x}(k) \\
            \end{array}
            \right]'\left[
            \begin{array}{cc}
                0 & S_1 \\
                S_1' & S_2 \\
            \end{array}
            \right]\left[
            \begin{array}{c}
                x(k) \\
                \tilde{x}(k) \\
            \end{array}
            \right].\nonumber\\
            \label{o21}
        \end{eqnarray}
        Moreover, the difference between the corresponding cost function (\ref{o21}) and the optimal cost under state feedback controller (all information are sharing to all controllers) is given by
        \begin{eqnarray}
            \delta J(s,\infty)&\triangleq& J^\star(s,\infty)-J^*(s,\infty)\nonumber\\
            &=&\sum_{k=s}^\infty\left[
            \begin{array}{c}
                x(k) \\
                \tilde{x}(k) \\
            \end{array}
            \right]'\left[
            \begin{array}{cc}
                0 & S_1 \\
                S_1' & S_2 \\
            \end{array}
            \right]\left[
            \begin{array}{c}
                x(k) \\
                \tilde{x}(k) \\
            \end{array}
            \right].\label{o23}
        \end{eqnarray}

    \end{theorem}
    \emph{Proof.}
    From (\ref{o13}) and (\ref{o8}), we have that
    \begin{eqnarray}
        &&x'(k)Px(k)-x'(k+1)Px(k+1)\nonumber\\
        &=&x'(k)[P-(A+BK)'P(A+BK)]x(k)\nonumber\\
        &&-\tilde{x}'(k)\mathcal{B}'P(A+BK)x(k)\nonumber\\
        &&-x'(k)(A+BK)'P\mathcal{B}\tilde{x}(k)\nonumber\\
        &&-\tilde{x}'(k)\mathcal{B}'P\mathcal{B}\tilde{x}(k)\nonumber\\
        &=&x'(k)(Q+K'RK)x(k)\nonumber\\
        &&-\tilde{x}'(k)\mathcal{B}'P(A+BK)x(k)\nonumber\\
        &&-x'(k)(A+BK)'P\mathcal{B}\tilde{x}(k)\nonumber\\
        &&-\tilde{x}'(k)\mathcal{B}'P\mathcal{B}\tilde{x}(k).\nonumber
    \end{eqnarray}
    By taking summation on $k$ from $s$ to $M$ and applying algebraic calculations, there holds that
    \begin{eqnarray}
        &&x'(s)Px(s)-x'(M+1)Px(M+1)\nonumber\\
        &=&J(s, M)-\sum_{k=s}^M\left[
        \begin{array}{c}
            x(k) \\
            \tilde{x}(k) \\
        \end{array}
        \right]'\left[
        \begin{array}{cc}
            0 & S_1 \\
            S_1' & S_2 \\
        \end{array}
        \right]\left[
        \begin{array}{c}
            x(k) \\
            \tilde{x}(k) \\
        \end{array}
        \right].\nonumber
    \end{eqnarray}
    Together with the fact that $\lim_{M\rightarrow\infty}x'(M+1)Px(M+1)=0$ obtained in Theorem  \ref{t1}, (\ref{o20}) follows
    by letting $M$ tend to $\infty$.

    Furthermore, from (\ref{o22}), the optimal cost under centralized controller
    is given by
    \begin{eqnarray}
        J^*(s,\infty)=x'(s)Px(s).\nonumber
    \end{eqnarray}
    Combining with (\ref{o21}), the difference (\ref{o23}) follows.
    The proof is now completed. \hfill $\blacksquare$

Finally, we analyze the asymptotical optimal property of the decentralized controllers (\ref{o11}) and (\ref{o12}).
 In fact, from (\ref{o23}), the derivation of the corresponding cost function under decentralized controllers (\ref{o11}) and (\ref{o12})
 from the optimal cost depends strongly on $\left[
 \begin{array}{c}
     x(k) \\
     \tilde{x}(k) \\
 \end{array}
 \right]$. Recalling the system (\ref{o19}) and the stability of the matrix $\bar{A}$, there exist positive constants
 $\lambda<1$ and $c$ such that
 \begin{eqnarray}
     \Big\|\left[
     \begin{array}{c}
         x(k) \\
         \tilde{x}(k) \\
     \end{array}
     \right]\Big\|\leq c\lambda^k\Big\|\left[
     \begin{array}{c}
         x(0) \\
         \tilde{x}(0) \\
     \end{array}
     \right]\Big\|.\label{o34}
 \end{eqnarray}
 Thus, we have
 \begin{eqnarray}
     &&\sum_{k=s}^\infty\left[
     \begin{array}{c}
         x(k) \\
         \tilde{x}(k) \\
     \end{array}
     \right]'\left[
     \begin{array}{cc}
         0 & S_1 \\
         S_1' & S_2 \\
     \end{array}
     \right]\left[
     \begin{array}{c}
         x(k) \\
         \tilde{x}(k) \\
     \end{array}
     \right]\nonumber\\
     &\leq &\sum_{k=s}^\infty\Big\|\left[
     \begin{array}{cc}
         0 & S_1 \\
         S_1' & S_2 \\
     \end{array}
     \right]\Big\|\Big\|\left[
     \begin{array}{c}
         x(k) \\
         \tilde{x}(k) \\
     \end{array}
     \right]\Big\|^2\nonumber\\
     &\leq &\Big\|\left[
     \begin{array}{cc}
         0 & S_1 \\
         S_1' & S_2 \\
     \end{array}
     \right]\Big\|\Big\|\left[
     \begin{array}{c}
         x(0) \\
         \tilde{x}(0) \\
     \end{array}
     \right]\Big\|^2c^2\sum_{k=s}^\infty \lambda^{2k}\nonumber\\
     &= &\bar{c}\lambda^{2s},\label{o24}
 \end{eqnarray}
 where $\bar{c}\triangleq\Big\|\left[
 \begin{array}{cc}
     0 & S_1 \\
     S_1' & S_2 \\
 \end{array}
 \right]\Big\|\Big\|\left[
 \begin{array}{c}
     x(0) \\
     \tilde{x}(0) \\
 \end{array}
 \right]\Big\|^2\frac{c^2}{1-\lambda^2}$.

Noting that $0<\lambda<1$, for any $\varepsilon>0$, there exists a sufficiently large integer $N$ such that
 \begin{eqnarray}
     \lambda^{2N}<\frac{\varepsilon}{\bar{c}}.\nonumber
 \end{eqnarray}
 Together with (\ref{o24}), it holds that
 \begin{eqnarray}
     \sum_{k=N}^\infty\left[
     \begin{array}{c}
         x(k) \\
         \tilde{x}(k) \\
     \end{array}
     \right]'\left[
     \begin{array}{cc}
         0 & S_1 \\
         S_1' & S_2 \\
     \end{array}
     \right]\left[
     \begin{array}{c}
         x(k) \\
         \tilde{x}(k) \\
     \end{array}
     \right]< \varepsilon.\label{o33}
 \end{eqnarray}
 In this sense, we say that the decentralized controllers (\ref{o11}) and (\ref{o12}) are asymptotically optimal, that is,
 there exists a sufficiently large integer $N$ such that the difference (\ref{o23}) between the corresponding cost function (\ref{o21}) and the optimal cost under state feedback controller satisfies
 \begin{eqnarray}
     \delta J(N,\infty)<\varepsilon.\nonumber
 \end{eqnarray}

    \section{Extension to the Case with Multiple Inputs}

    In this section, we consider the system with multiple inputs as
    \begin{eqnarray}
        x(k+1)&=&Ax(k)+\sum_{i=1}^rB_iu_i(k),\label{o34}\\
        y_i(k)&=&H_ix(k), i=1,2, \ldots, r,\label{o35}
    \end{eqnarray}
    and the cost function is given by
    \begin{eqnarray}
        J&=&\sum_{k=0}^\infty[x'(k)Qx(k)+u_1'(k)R_1u_1(k)+\ldots\nonumber\\
        &&~~~~~~~+u_r'(k)R_ru_r(k)],\label{o36}
    \end{eqnarray}
where $B_i\in R^{n\times m_i}$, $H_i\in R^{s_i\times n}$ are constant matrices, and $R_1,\ldots,R_r$ are positive definite matrices with appropriate dimensions.

    Similar to (\ref{o3}) and (\ref{o25}), we define for $i=1,\ldots, r$ that
    \begin{eqnarray}
        \bar{Y}_i(k)&=&\{y_i(0), \ldots, y_i(k), u_1(0), \ldots, u_1(k-1), \nonumber\\
        &&~~~~~~~~~~~~\ldots, u_r(0), \ldots, u_r(k-1)\}, \label{o38}\\
        Y_i(k)&=&\{y_i(0), \ldots, y_i(k), u_i(0), \ldots, u_i(k-1)\}. \label{o37}
    \end{eqnarray}
    By applying same procedures to the derivations in Section III-V, the solutions to the scenarios with the $\bar{Y}_i(k)$-causal controller and the $Y_i(k)$-causal controller are given below.

    Firstly, we give the solution of the $\bar{Y}_i(k)$-causal controller (i.e., historical controllers are sharing) for $i=1,\ldots, r$
   by defining the decentralized controllers :
    \begin{eqnarray}
        u_i^\diamond(k)=\mathcal{K}_i\hat{x}_i^\diamond(k),\label{o40}
    \end{eqnarray}
    and the observers $\hat{x}_i^\diamond$:
    \begin{eqnarray}
        \hat{x}_i^\diamond(k+1)&=&A\hat{x}_i^\diamond(k)+B_1u_1^\diamond(k)+\ldots+B_ru_r^\diamond(k)\nonumber\\
        &&+L_i^\diamond[y_i(k)-H_i\hat{x}_i^\diamond(k)],\label{o39}
    \end{eqnarray}
    In particular, the feedback gains are given by  $\mathcal{K}_1=\left[
    \begin{array}{cccc}
        I & 0 &\cdots &0\\
    \end{array}
    \right]\mathcal{K},
    \cdots, \mathcal{K}_r=\left[
    \begin{array}{cccc}
        0&\ldots &0 & I \\
    \end{array}
    \right]\mathcal{K}$
    where $\mathcal{K}$ is defined by
    \begin{eqnarray}
        \mathcal{K}&=&-(\mathcal{R}+\mathbb{B}'\mathcal{P}\mathbb{B})^{-1}\mathbb{B}'\mathcal{P}A,\label{o41}
    \end{eqnarray}
    and $\mathcal{P}$ satisfies ARE:
    \begin{eqnarray}
        \mathcal{P}&=&A'\mathcal{P}A+Q-A'\mathcal{P}\mathbb{B}(\mathcal{R}+\mathbb{B}'\mathcal{P}\mathbb{B})^{-1}\mathbb{B}'\mathcal{P}A,\label{o42}
    \end{eqnarray}
    with $\mathbb{B}=\left[
    \begin{array}{ccc}
        B_1 & \cdots & B_r \\
    \end{array}
    \right]
    $ and $\mathcal{R}=diag\{R_1, R_2,\ldots, R_r\}.$

    \begin{theorem}\label{t5}
        Assume that $(A,\mathbb{B})$ is stabilizable, systems $(A,Q)$ and $(A, H_i)$ are observable for $i=1, \ldots, r$, then the closed-loop system (\ref{o34}) under controller (\ref{o40}) is stable:
        \begin{eqnarray}
            x(k+1)&=&Ax(k)+\sum_{i=1}^rB_i\mathcal{K}_i\hat{x}_i^\diamond(k).\nonumber
        \end{eqnarray}

    \end{theorem}
    \emph{Proof.} The proof is similar to that in Theorem \ref{t5}. So we omit the details. \hfill $\blacksquare$

    Secondly, we give the solution of the $Y_i(k)$-causal controller (strictly private information)
    for $i=1,\ldots, r$ by defining the decentralized controllers:
    \begin{eqnarray}
        u_i^\star(k)=\mathcal{K}_i\hat{x}_i(k),\label{o43}
    \end{eqnarray}
    where the observers $\hat{x}_i$ satisfy
    \begin{eqnarray}
        \hat{x}_1(k+1)&=&A\hat{x}_1(k)+B_1u_1^\star(k)+B_2\mathcal{K}_2\hat{x}_1(k)+\cdots\nonumber\\
        &&+B_r\mathcal{K}_r\hat{x}_1(k)+\mathcal{L}_1[y_1(k)-H_1\hat{x}_1(k)],\label{o44}\\
        &&\vdots\nonumber\\
        \hat{x}_r(k+1)&=&A\hat{x}_r(k)+B_1\mathcal{K}_1\hat{x}_r(k)+\cdots\nonumber\\
        &&+B_{r-1}\mathcal{K}_{r-1}\hat{x}_r(k)+B_ru_r^\star(k)\nonumber\\
        &&+\mathcal{L}_r[y_r(k)-H_r\hat{x}_r(k)].\label{o45}
    \end{eqnarray}
By defining the matrix $\mathbb{A}$ as given by (\ref{o440}),
\begin{figure*}[ht]
\hrulefill
 \begin{equation}
    \mathbb{A}=\left[
    \begin{array}{cccc}
        A+\mathbb{B}K-B_1K_1-L_1H_1 & -B_2K_2 & \cdots & -B_rK_r\\
        -B_1K_1 & A+\mathbb{B}K-B_2K_2-L_2H_2 & \cdots & -B_rK_r \\
        &  & \cdots &  \\
        -B_1K_1 & \cdots & -B_{r-1}K_{r-1} & A+\mathbb{B}K-B_rK_r-L_rH_r \\
    \end{array}
    \right]\label{o440}
\end{equation}
\end{figure*} the solution of the private information case is shown below.
    \begin{theorem}
        If the observer gains $\mathcal{L}_i, i=1,\ldots, r$ are chosen such that the matrix $\mathbb{A}$ is stable,
        then the observers $\hat{x}_i(k), i=1,\ldots, r$ in (\ref{o44})-(\ref{o45}) are stable under controllers (\ref{o43}).
        Moreover, under the stabilizability of system $(A,\mathbb{B})$ and the observability of system $(A,Q)$,
        the closed-loop system (\ref{o34}) under controllers (\ref{o43}) is stable. The associated performance reaches asymptotically optimal.
    \end{theorem}
    \emph{Proof.} The proof is similar to that in Theorem \ref{t1}. So we omit the details. \hfill $\blacksquare$

    \section{Numerical Examples}

    In this section, we present an example to verify the effectiveness of the derived results.
    Consider system (\ref{o1}) with parameters given by
    \begin{eqnarray}
        &&A=\left[
        \begin{array}{cc}
            1 & 1 \\
            2 & -1 \\
        \end{array}
        \right], B_1=\left[
        \begin{array}{c}
            0.6 \\
            0.5 \\
        \end{array}
        \right],B_1=\left[
        \begin{array}{c}
            0 \\
            1 \\
        \end{array}
        \right],\nonumber\\
        &&H_1=\left[
        \begin{array}{cc}
            1 & 0 \\
        \end{array}
        \right], H_2=\left[
        \begin{array}{cc}
            0 & 1 \\
        \end{array}
        \right], \nonumber\\
        &&Q=I,R_1=1, R_2=1.\nonumber
    \end{eqnarray}
    From solving ARE (\ref{o8}), we have the feedback gain in (\ref{o7}) as
    \begin{eqnarray}
        K=\left[
        \begin{array}{cc}
          -1.2382 &  -0.7982\\
        -1.1262   & 1.1412\\
        \end{array}
        \right].\nonumber
    \end{eqnarray}
    By applying Theorem \ref{t1}, we select $L_1=\left[
    \begin{array}{cc}
        0.3 & 0.5\\
    \end{array}
    \right]
    , L_2=\left[
    \begin{array}{cc}
        0.8 & -0.6\\
    \end{array}
    \right]$,
    the errors (\ref{o16}) between the state and the observers (\ref{o9}) and (\ref{o10}) are stable as shown in Fig. 1.
    Moreover, the closed-loop system (\ref{o1}) under controllers (\ref{o11})-(\ref{o12}) is stable
    as shown in Fig. 2. It is seen from Fig. 1 and 2 that $N$ in (\ref{o33}) can be chosen as $20$ to achieve the asymptotical optimality.

    \begin{figure}
        \centering
        \includegraphics[width=70mm]{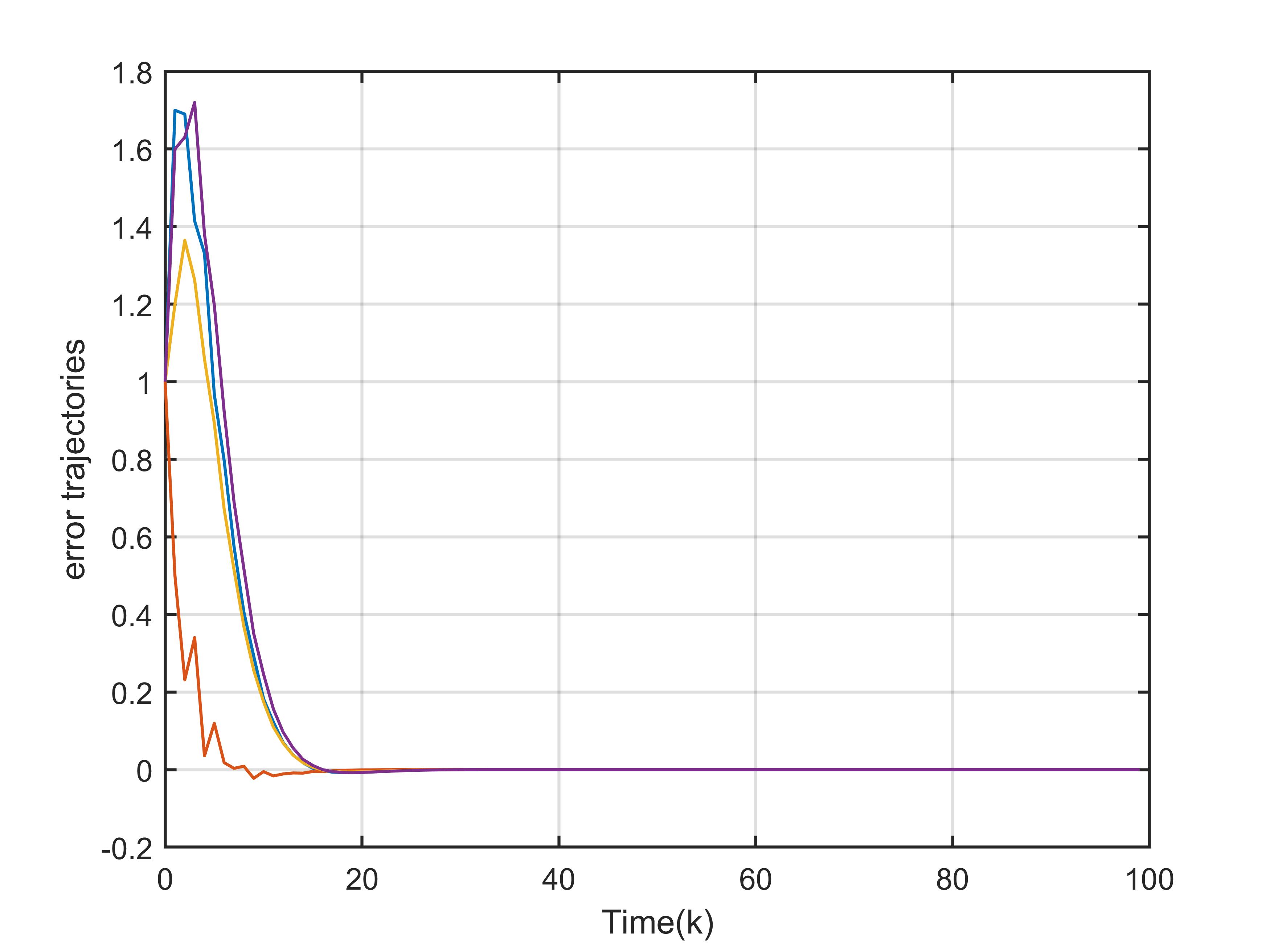}\\
        \caption{The trajectory of $\tilde{x}(k)$.}
    \end{figure}

    \begin{figure}
        \centering
        \includegraphics[width=70mm]{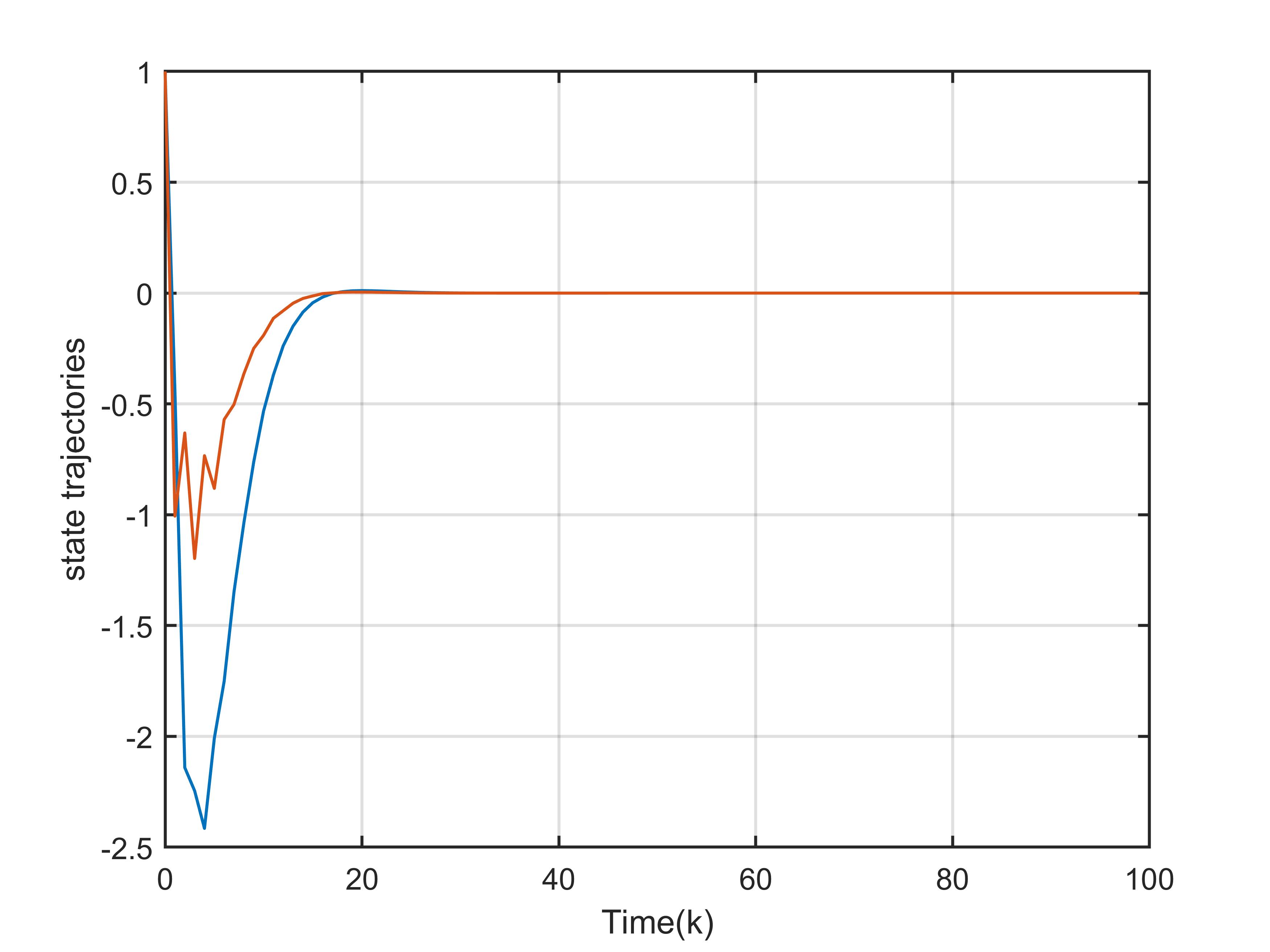}\\
        \caption{The trajectory of $x(k)$.}
    \end{figure}

    \section{Conclusions}

  This paper was concerned with the LQ optimal control problem of linear systems with private input and measurement information.
    To overcome the difficulty from the unavailability of other regulators' historical input information, we defined
    a novel kind of observers and decentralized controllers by using only local information, which was proven that the corresponding cost function under the proposed decentralized controllers are asymptotically optimal as comparison with the optimal cost under optimal controller with full state information.
    Numerical examples verified the effectiveness of the derived results.


\end{document}